\documentstyle[11pt]{amsart}

\title[Manifolds of positive isotropic curvature]{The fundamental group of manifolds of positive isotropic curvature and surface groups}
\author{Ailana Fraser}
\address{Department of Mathematics \\
                 University of British Columbia \\
                 Vancouver, BC V6T 1Z2}
\author{Jon Wolfson}
\address{Department of Mathematics \\ 
                 Michigan State University \\
                 East Lansing, MI 48824}

\thanks{  The first author was partially supported by  the
Natural Sciences and Engineering Research Council of Canada (NSERC) and the second 
author was partially supported by NSF grant DMS-0304587 }

\email{afraser@@math.ubc.ca and wolfson@@math.msu.edu}

\date{January 10, 2005}

\newtheorem{thm}{Theorem}[section]
\newtheorem{lem}[thm]{Lemma}

\theoremstyle{definition}
\newtheorem{rem}{Remark}[section]

\newtheorem{defn}{Definition}[section]

\numberwithin{equation}{section}

\renewcommand{\a}{\alpha}

\newcommand{\e}{\varepsilon}
\newcommand{\g}{\gamma}
\newcommand{\G}{\Gamma}

\newcommand{\n}{\nabla}
\newcommand{\p}{\partial}

\newcommand{\Sig}{\Sigma}

\renewcommand{\t}{\tau}

\def\Pb{\ifmmode{\mathbb P}\else{$\mathbb P$}\fi}
\def\Z{\ifmmode{\mathbb Z}\else{$\mathbb Z$}\fi}
\def\Q{\ifmmode{\mathbb Q}\else{$\mathbb Q$}\fi}
\def\C{\ifmmode{\mathbb C}\else{$\mathbb C$}\fi}
\def\R{\ifmmode{\mathbb R}\else{$\mathbb R$}\fi}
\def\S{\ifmmode{S^2}\else{$S^2$}\fi}

\newcommand{\la}{\langle}
\newcommand{\ra}{\rangle}
\newcommand{\der}[1]{\frac{\partial}{\partial #1}}
\newcommand{\fder}[2]{\frac{\partial #1}{\partial #2}}

\def\Ca{\cal C}

\def\S{\cal S}

\begin{document}

\maketitle

\begin{abstract}
In this paper we study the topology of  compact manifolds of positive 
isotropic curvature (PIC). There are many examples of  non-simply 
connected compact manifolds with  positive isotropic curvature. We 
prove that the fundamental group of a compact Riemannian manifold 
with PIC, of dimension $\geq 5$, does not contain  a subgroup 
isomorphic to the fundamental group of a compact Riemann surface. The 
proof uses stable minimal surface theory.
\end{abstract}

\setcounter{secnumdepth}{1}

\setcounter{section}{-1}

\section{Introduction}

In this paper we study fundamental groups of compact manifolds of positive isotropic curvature.
We prove that the fundamental group of a 
compact Riemannian manifold with positive isotropic curvature of dimension $\geq 5$ can not contain a surface group as a
subgroup:

\begin{thm} \label{thm:pic2}
Let $M$ be a compact $n$-dimensional Riemannian manifold, $n \geq 5$, with positive
isotropic curvature. Then the fundamental group of $M$, $\pi_1(M)$, does not contain a subgroup
isomorphic to the fundamental group $\pi_1(\Sig_0)$ of a compact Riemann 
surface $\Sig_0$ of genus $g_0  \geq 1$.
\end{thm}

In ~\cite{F} the first author proved the genus one case. In particular she proved that the the fundamental group of 
a compact $n$-manifold, $n \geq 5$, of positive isotropic curvature does not contain a subgroup isomorphic
to $\Z \oplus \Z$. The proof we give here of Theorem \ref{thm:pic2} is closely modeled
on the proof in~\cite{F}.

The main open conjecture on the topology of compact  Riemannian manifolds with positive
isotropic curvature concerns the fundamental group.
It is conjectured that the fundamental group of  a compact  Riemannian manifold with positive
isotropic curvature is virtually free (i.e., contains a free subgroup of finite index). Recall that the 
fundamental group of  a compact  manifold has $0,1,2$ or infinitely many ends [E]. Here we can define
the number of ends of a fundamental group as the number of geometric ends of the universal cover.
The conjecture implies that virtually (i.e., up to finite covers) every subgroup of the fundamental group 
of  a compact  Riemannian manifold with PIC has either $2$ or infinitely many ends. On the other hand, the 
fundamental group of a compact Riemann surface of genus $\geq 1$ has one end. In this sense the theorem can be
seen as evidence for this conjecture. Along this line of reasoning we believe the following weaker conjecture
remains interesting and is more amenable: The fundamental group of  a compact  Riemannian manifold with positive
isotropic curvature has virtually no subgroup with exactly one end.

\bigskip

Recall the definition of positive isotropic curvature. Let $M$ be an $n$-dimensional Riemannian manifold. The inner product on the tangent space $T_pM$ at a point $p \in M$ can be extended to the complexified tangent space $T_pM \otimes \mathbb{C}$ as a complex bilinear form 
$( \cdot , \cdot )$ or as a Hermitian inner product $\la \cdot , \cdot \ra$. The relationship between 
these extensions is given by $\la v,w \ra = (v , \bar{w})$ for $v$, $w \in T_pM \otimes \mathbb{C}$.
The curvature tensor extends to complex vectors by linearity, and the {\em complex sectional curvature} of a two-dimensional subspace $\pi$ of $T_pM \otimes \mathbb{C}$ is defined by
$K(\pi)=\la R(v,w)\bar{w}, v \ra$, where $\{v , w \}$ is any unitary basis of $\pi$. A subspace 
$\pi \subset T_pM \otimes \mathbb{C}$ is said to be {\em isotropic} if every vector $v \in \pi$ has 
square zero; that is, $(v,v)=0$. 
\begin{defn}
A Riemannian manifold $M$ has {\em positive isotropic curvature}
(PIC) if $K(\pi) > 0$ for every isotropic two-plane $\pi \subset T_pM \otimes \mathbb{C}$,
for all $p \in M$.
\end{defn}
This curvature condition is nonvacuous only for $n \geq 4$, since in dimensions less than four 
there are no two-dimensional isotropic subspaces. PIC is a curvature condition that arises very
naturally when studying stability of minimal surfaces, just as positive sectional curvature is 
ideally adapted to studying stability of geodesics. Any manifold with pointwise quarter-pinched
sectional curvatures or positive curvature operator has PIC. PIC implies positive scalar curvature,
but not positive (or even nonnegative) Ricci curvature. For more background on this curvature
condition we refer, for example, to the introductions of ~\cite{M-M}, ~\cite{M-W} and ~\cite{F}.

Theorem \ref{thm:pic2} is proved by assuming the existence of a subgroup $G \subset \pi_1(M)$,
isomorphic to a surface group and deriving a contradiction. In an analogous situation Schoen and Yau ~\cite{S-Y}
derive a contradiction to the existence of such a subgroup of the fundamental group of a compact three manifold $X$
of positive scalar curvature. They construct a stable minimal surface $u: \Sigma_0 \to X$ with $u_* : \pi_1(\Sigma_0) \to  \pi_1(X)$
a monomorphism and use the curvature condition to derive a contradiction to the stability of $u$. If $L$ is a lens space
then the manifold $S^1 \times L$ admits a PIC metric and a stable minimal map $u: T^2 \to S^1 \times L$. Therefore
an argument like that of ~\cite{S-Y} cannot be expected to work for manifolds with PIC. Rather we assume the existence of a  subgroup
$G$ isomorphic to a surface group and find a suitable finite index normal subgroup $N$ of $G$ using covering space theory. 
A contradiction results from the existence of a stable minimal
map $u: \Sigma \to M$ with $u_*:  \pi_1(\Sigma) \to  \pi_1(M)$ an isomorphism onto $N$.

\section{Proof of the theorem}


Let $M$ be a compact Riemannian manifold. If $\a \in \pi_1(M)$ we define the {\it systole of $\a$} to be:
$${\cal{S}}(\a) = \inf \{ \ell(\g): \g \;\; \mbox{is a closed rectifiable curve with} \;\; [\g]=\a \}.$$
where $\ell(\g)$ denotes the length of $\g$. We define the {\it systole of $M$} to be:
$${\cal{S}}(M)= \inf \{ {\cal{S}}(\a) : \a \in \pi_1(M), \;\;  \a \neq 0 \}.$$

\begin{thm}
\label{thm:cover}
Let $M$ be a  compact Riemannian manifold whose fundamental group $\pi_1(M)$ has a subgroup $G$ isomorphic to $\pi_1(\Sig_0)$, where $\Sig_0$ is a compact Riemann surface of genus $g_0 \geq 2$. Given $C > 0$, there is an integer $k > 0$ and a  index $k$ normal subgroup $N$ of $G$ such that: (i)  there is a smooth map  $h: \Sig \to M$ of a compact surface into $M$ with $h_*:\pi_1(\Sig) \to \pi_1(M)$ a monomorphism onto $N$ , (ii) for any such map $h$, with respect to the induced metric, the systole of $\Sig$ is $>C$ (i.e., every closed non-trivial geodesic $\g$ has length $\ell(\g) > C$).
\end{thm}

\begin{pf}
Given $C > 0$ there are at most finitely many free homotopy classes $\{ \g_i \in G \subset  \pi_1(M): i=1,\dots,k \}$ with systole $\leq C$. The fundamental group of a compact Riemann
surface is residually finite~\cite{M-K-S}. Therefore for each class $[\g_i]$ there is a finite index normal subgroup $N_i$ of $G$ such that $[\g_i] \notin N_i$. The intersection $N = \cap_{i=1}^k N_i$ is a finite index (index $k$) normal subgroup of $G$. Let $h_0: \Sig_0 \to M$ be a smooth map such that ${h_0}_*:\pi_1(\Sig_0) \to G$ is an isomorphism. There is a regular $k$-covering $p: \Sig \to \Sig_0$ such that $p_*: \pi_1(\Sig) \to N$ is an isomorphism. The map $h =  h_0  \circ p$ satisfies (i) in the statement of the theorem. For any map $h$ with $h_*: \pi_1(\Sig) \to N$
and any non-trivial closed curve $\g$, $h_*([\g]) \in N$. Therefore $\ell(\g) > C$, where the length is computed using the induced metric. The result follows. 
 \end{pf}

\begin{thm}
\label{thm:distancemap}
Given $\e > 0$ there exists $C(\e) > 0$ such that if the systole of the compact Riemann surface $\Sig$ 
is $> C(\e)$ then there is a  Lipschitz distance decreasing, degree one map $f: \Sig \to S^2$ satisfying $|df| < \e$.
\end{thm}

This can be proved as in~\cite{F} in the proof of Lemma 3.1. For completeness, we repeat the
argument here.

\begin{pf}
Let $\tilde{\Sig}$ be the universal cover of $\Sig$, and
$p:\tilde{\Sig} \rightarrow \Sig$ the covering map.
Since $\tilde{\Sig}$ is complete and noncompact with compact quotient
$\Sig$, there is a geodesic line 
$r: \R \rightarrow \tilde{\Sig}$.
Let $\Ca$ be the component of $\tilde{\Sig}-r(\R)$
which is on the side of $\tilde{\Sig}$ in the direction of the unit 
normal $\nu$ to $r$ such that $\{r'(0), \nu(0) \}$ is positively oriented.
Choose $T$ very large, $T>>C$.
Define $D_1: \tilde{\Sig} \rightarrow \R$ by,
\[
     D_1(x)=d(x,r(T))-T
\]
and let $D_2: \tilde{\Sig} \rightarrow \R$ be the signed
distance function to $r$,
\[
   D_2(x)= \begin{cases}
    d(x,r)& \text{for}\;x \in \Ca   \\
    -d(x,r)& \text{for}\;x \in \tilde{\Sig} - \Ca
    \end{cases}
\]
Both $D_1$ and $D_2$ are Lipschitz continuous with derivative
bounded by 1. 
Consider the region
\[
   \tilde{\mathcal{R}}=\Big\{x \in \tilde{\Sig} : 
   |D_1(x)| \leq \frac{C}{4},
   |D_2(x)| \leq \frac{C}{4} \Big\}
\]
Define the map $F: \tilde{\mathcal{R}} \rightarrow \R^2$ by
\[
     F(x)= (D_1(x),D_2(x)).
\]
Then the boundary of $\tilde{\mathcal{R}}$ is mapped by $F$ to the boundary of the 
rectangle
$[-\frac{C}{4},\frac{C}{4}]\times[-\frac{C}{4},\frac{C}{4}]$ 
in $\R^2$.
Also, $r(0)$ is the only point in $\tilde{\mathcal{R}}$ which is mapped under 
$F$ to the origin in $\R^2$. 
In fact, $F$ is a local diffeomorphism in a neighborhood of $r(0)$, 
and hence the degree of $F$ is equal to one on the component of 
$\R^2-F(\p \tilde{\mathcal{R}})$ containing the origin.
Hence, $F(\tilde{\mathcal{R}})$ covers a disk of radius at least 
$\frac{C}{4}$ about the origin in $\R^2$.

Let $\lambda : D(0,\frac{C}{4}) \rightarrow D(0,\pi)$ be the
contraction $\lambda(x)=\frac{4\pi}{C}x$ (where $D(0,s)$ denotes
the disk of radius $s$ centered at the origin in $\R^2$). Let 
$e: D(0,\pi) \rightarrow S^2$ be the exponential map at the north
pole $n$ of the sphere $S^2$ with the standard metric, 
$e(x)=\exp_n(x)$. Extend $g=e \circ \lambda$ to $F(\tilde{\mathcal{R}})$
by defining $g \equiv s$ (the south pole) on 
$F(\tilde{\mathcal{R}}) - D(0,\frac{C}{4})$. Then 
$\tilde{f}=g \circ F : \tilde{\mathcal{R}} \rightarrow S^2$ 
is Lipschitz with derivative
\[
     |d\tilde{f}| \leq |dF||d\lambda||de| \leq \frac{c_1}{C}
\]
where $c_1$ is a constant (independent of $\Sig$).

Observe that 
diam $\tilde{\mathcal{R}} \leq C$. Given $x, y \in \tilde{\mathcal{R}}$,
we have
$d(x,r)=d(x,r(t_1)) \leq \frac{C}{4}$ and
$d(y,r)=d(y,r(t_2)) \leq \frac{C}{4}$ 
for some $t_1, t_2$ with $r(t_1), r(t_2) \in \tilde{\mathcal{R}}$, 
and
\[
    d(r(t_1),r(t_2)) \leq D_1(r(t_1))
                                 +  D_1(r(t_2))
    \leq \frac{C}{2}.
\]
Now,
\[
    d(x,y) \leq d(x, r(t_1)) + d(r(t_1),r(t_2))
                 + d(r(t_2),y) \leq C.
\]
Hence no two points of $\tilde{\mathcal{R}}$ are identified in the quotient
$\Sig$. If two points were identified, then any minimal curve in
$\tilde{\mathcal{R}}$ joining the two points would project to a nontrivial 
closed curve of length less than or equal to $C$, which is less 
than the systole, a contradiction.

Therefore, $\tilde{\mathcal{R}}$ projects one to one into $\Sig$, and
we may define $f:   {\mathcal{R}}   \rightarrow S^2$, 
where $ {\mathcal{R}} =p(\tilde{\mathcal{R}})$, by
$f(x)=\tilde{f}(p^{-1}(x))$. 
Since $f$ maps $\p \mathcal{R}$ to $s$, we can extend  
$f$ from $\mathcal{R}$ to a map $f: \Sig \rightarrow S^2$ by defining 
$f \equiv s$ on $\Sig - \mathcal{R}$. \end{pf}

\noindent {\bf Remark:} Gromov-Lawson~\cite{G-L} use a similar idea to construct what they call {\it a c-contracting map} (\cite{G-L}, Proposition 3.3): Assume that a compact Riemannian manifold $(X, g)$ has residually finite fundamental group and admits a metric ${\tilde g}$ of non-positive curvature. Their result constructs a contracting map $X' \to S^n$, where $X'$ is some finite cover of $X$. The argument they use applies to the fixed metric $g$ and not (at least directly) to a family of metrics, as in our case.

\bigskip

Let $M$ be a Riemannian manifold of dimension $n \geq 5$ with positive isotropic curvature $> \kappa$. Suppose that the fundamental group $\pi_1(M)$ has a subgroup $G$ isomorphic to $\pi_1(\Sig_0)$, where $\Sig_0$ is a compact Riemann surface of genus $g_0 \geq 2$. Let $h: \Sig_0 \to M$ be a smooth map such that 
${h_0}_*:\pi_1(\Sig_0) \to \pi_1(M)$ is an isomorphism onto $G$.  
Choose $\e$ such that $0 < 3(c\varepsilon)^2 < \kappa$ where
$c$ is to be chosen later. Let $C(\e)$ be given by Theorem \ref{thm:distancemap}.
Using Theorem \ref{thm:cover} there is a compact Riemann surface $\Sig$ of genus $g$,  a regular $k$-covering $p: \Sig \to \Sig_0$ and a smooth map $h: \Sig \to M$ whose induced map on $\pi_1$ is injective. 
Since $p$ is a covering map,
the Euler characteristics satisfy $\chi(\Sig) = k \chi(\Sig_0)$ and therefore $2 - 2g = k(2-2g_0)$. Also by Theorem \ref{thm:cover} for every map $\tilde{h}: \Sig \to M$ whose induced map on $\pi_1(\Sig)$ equals $h_*$ the surface $\Sig$, with the induced metric, has systole $> C(\e)$.
The map $h$ is incompressible and thus following~\cite{S-Y} there exists a stable conformal branched minimal immersion $u: \Sig \to M$ whose induced map on  $\pi_1(\Sig)$ is $h_*$.

Let  $u^*({\mathcal{N}})$ be the pull-back of the 
normal bundle of the minimal surface, $u(\Sig)$, with the pull back of the metric and
normal connection $\nabla^{\perp}$.
Let $E=u^*({\mathcal{N}}) \otimes \mathbb{C}$ be the complexified normal bundle.
Then $c_1(E) = 0$ since $E$ is the complexification of a real bundle.
The metric on $u^*({\mathcal{N}})$ extends as a complex bilinear form
$( \cdot , \cdot )$ or as a Hermitian metric $\la \cdot , \cdot \ra$
on $E$, and the connection $\n^{\perp}$ and curvature tensor extend
complex linearly to sections of $E$.
There is a unique holomorphic structure on $E$ such that the
$\bar{\p}$ operator is given by
      $\bar{\p}\omega=(\nabla_{\der{\bar{z}}}^{\perp} \omega)d\bar{z}$
where $\der{\bar{z}}=\frac{1}{2}(\der{x}+i\der{y})$, in local coordinates
$x, y$ on $\Sig$.
Using the complexified
formula for the second variation of area (see \cite{Si-Y},  \cite{M-M}, \cite{F})
the stability condition is the inequality:
\[
    \int_{\Sig} \la R(s,\fder{u}{{z}})\fder{u}{\bar{z}},s \ra \; dxdy
    \leq \int_{\Sig} [|\nabla_{\der{\bar{z}}}^{\perp} s|^2 
    - |\nabla_{\der{z}}^{\top} s|^2] \; dxdy
\]
for all $s \in \Gamma(E)$. Assume that $s$ is isotropic. Since $u$
is conformal, $\fder{u}{z}$ is isotropic and $\{s,\fder{u}{z}\}$ span
an isotropic two-plane. Using the lower bound on the isotropic
curvature and throwing away the second term on the right, we have
\begin{equation} \label{equation:stability}
    \kappa \int_{\Sig} |s|^2 \; da 
    \leq \int_{\Sig} |\bar{\p} s|^2 \; da
\end{equation}
where $da$ denotes area element for the induced metric $u^*g$ on $\Sig$.
We now argue that we can find an ``almost holomorphic'' isotropic
section of $E$ that violates this stability inequality.
That is, we will find $s \in \Gamma(E)$ such that
\[
     \int_{\Sig} |\bar{\p} s|^2 \; da 
    < \kappa \int_{\Sig} |s|^2 \; da
\]
The contradiction (to the existence of $u: \Sig \to M$) proves Theorem \ref{thm:pic2}.

By Theorem \ref{thm:distancemap}, since the systole of $\Sig$ is $>C(\e)$,  there is a distance decreasing, degree one map
$f: \Sig \rightarrow S^2$ with $|df| < \e$.
Let $(L_k, \nabla)$ be a complex line bundle over $S^2$ with metric and connection, 
with $c_1(L_k) > 3k(g_0-1)+1$.
Let $\xi = f^*L_k$ be the pulled back bundle over $\Sig$.
We denote the induced connection $\nabla$ and note that it defines a holomorphic structure on $\xi$,
and $c_1(\xi) > 3k(g_0-1)+1$.
The tensor product bundle $\xi \otimes E$ is a holomorphic rank $(n-2)$ bundle
over $\Sig$. Let ${\cal{H}} (\xi \otimes E)$ denote the complex vector space of holomorphic sections of $\xi \otimes E$. By the Riemann-Roch theorem:
\begin{equation}
\begin{split}
   {\rm{dim}}_{\C}\; {\cal{H}} (\xi \otimes E) & \geq  c_1( \mbox{det} (\xi \otimes E)) + (n-2) (1-g) \\
   & = (n-2) c_1(\xi) + c_1(E) + (n-2) k (1-g_0) \\
   & = (n-2) [ c_1(\xi) + k (1- g_0) ] 
\end{split}
\label{equation:dim}
\end{equation}

As in \cite{F}, we may prove existence of a holomorphic and {\it isotropic} section of $\xi \otimes E$. 
Set  $m = \mbox{dim}_{\C}\;  {\cal{H}} (\xi \otimes E)$ and observe that the complex bilinear pairing on $E$ defines a complex bilinear pairing
\[
    {\mathcal{H}} (\xi \otimes E) \times {\mathcal{H}}(\xi \otimes E)
    \rightarrow {\mathcal{H}} (\xi \otimes \xi).
\]
Given any $x \in \Sig$ we obtain a homogeneous polynomial on 
${\mathbb{C}}^m \cong {\mathcal{H}}(\xi \otimes E)$ given by
$P_x(\sigma)=(\sigma, \sigma)(x)$.
The zero set
\[
   V(P_x)=\{ \sigma \in {\mathbb{P}}^{m-1}: P_x(\sigma)=0\}
\]
is a hypersurface in ${\mathbb{P}}^{m-1}$.
Now given $m-1$ distinct points, we obtain $m-1$ hypersurfaces
in ${\mathbb{P}}^{m-1}$, and observe that $m-1$ such hypersurfaces
in ${\mathbb{P}}^{m-1}$ intersect in a nonempty set of points.
The intersection is a set of holomorphic sections of 
$\xi \otimes E$ which are isotropic at $m-1$ distinct points.
Let $\sigma \in \mathcal{H}(\xi \otimes E)$ be such a section. 
Then $(\sigma, \sigma)$ is a holomorphic section of 
$\xi \otimes \xi$ with at least $m-1$ zeros. But the number
of zeros of a holomorphic section of $\xi \otimes \xi$ is
$2c_1(\xi)$. From (\ref{equation:dim})
\[
      m-1 \geq (n-2) [ c_1(\xi) + k (1- g_0) ] -1 > 2 c_1(\xi)
\]
since $n \geq 5$ and we chose the line bundle $L_k$ so that $c_1(L_k) \geq 3k(g_0-1)+1$.
It follows that $(\sigma, \sigma) \equiv 0$ and so $\sigma$ is
isotropic.
For brevity of notation we will suppress the $k$ in $L_k$.

We claim that any holomorphic isotropic section of $\xi \otimes E$ produces an 
{\em almost holomorphic} isotropic section of $E$.  This is achieved by applying an
almost parallel section of $\xi^*$, the pull back under the distance decreasing map $f$
of a section of $L^*$, to `undo' the line bundle part of the section. The argument of \cite{F}
does not require fine control of these dual sections. 
In the higher genus case, due to the dependence of $L$ on $k$, the construction of the
required dual sections is much more delicate, and uses the following result.

\begin{lem}
\label{lem:sections}
Given a line bundle $L \to S^2$ and a Borel measure $\mu$ on $S^2$
there are sections $t_1$, $t_2$ and a constant $c > 0$ (independent of $c_1(L)$ and $\mu$) such that:
\begin{enumerate}
\item $|t_1|^2 + |t_2|^2 \geq 1$. \\
\item $| \nabla t_1| \leq c$. \\
\item $\int_{S^2} |t_1|^2 d\mu \geq  \frac{1}{2} \int_{S^2} |t_2|^2 d\mu$.
\end{enumerate}
for some connection $\nabla$ on $L$. 
\end{lem}

\medskip

\begin{pf}
Fix  $r \ll 1$ independent of $c_1(L)$ and $\mu$. 
Consider the function $\mu: S^2 \to {\R_+}$ given by $\mu(p) = \mu(B_r(p))$. Let $p_1 \in S^2$ be a point where the minimum of $\mu$ is achieved. Without loss of generality we can assume that the line bundle $(L, \nabla)$ is flat on $S^2 \setminus B_{\frac{r}{2}}(p_1)$. Let $\t_1$ be a smooth section of $L$ that is parallel of length $\sqrt{2}$ on $S^2 \setminus B_{\frac{r}{2}}(p_1)$. Let $\t_2$ be a smooth section of $L$ such that $|\t_2|^2 < 2$ everywhere and $|\t_2| \geq 1$ on $B_r(p_1)$. Denote by $\phi_{r, p_1}$ a cut off function that vanishes in $B_{\frac{r}{2}}(p_1)$, is identically one on $S^2 \setminus B_r(p_1)$ and satisfies $| d\phi_{r, p_1}| \leq \frac{1}{r}$. Define:
$$t_1 = \phi_{r, p_1} \t_1, \;\;\; t_2 =  \t_2.$$
It follows that:
$$\int_{S^2 \setminus B_r(p_1)} |t_1|^2 d\mu \geq \int_{S^2 \setminus B_r(p_1)} |t_2|^2 d\mu.$$
By the choice of $p_1$, there is a point $p_2 \in S^2$ with $B_r(p_1)$ and $B_r(p_2)$ disjoint and $\mu(B_r(p_2)) \geq \mu(B_r(p_1))$. Then
$$\int_{B_r(p_2)} |t_1|^2 d\mu = 2\mu(B_r(p_2)) \geq 2\mu(B_r(p_1)) \geq \int_{B_r(p_1)} |t_2|^2 d\mu.$$
Therefore,
$$2 \int_{S^2} |t_1|^2 d\mu \geq  \int_{S^2} |t_2|^2 d\mu.$$
The result follows.

\end{pf}

Suppose that $\tilde{s} \in \G( \xi \otimes E)$ is holomorphic and isotropic.
Define the measure $\nu(A)= \int_A |\tilde{s}|^2 da$ for any $da$ measurable subset $A \subset  \Sig$.
Define the push forward measure:
\begin{equation} 
\mu = f_{\sharp} \nu,
\end{equation} 
Then $\mu$ is a Borel measure on $S^2$ such that for any function $h$ on $S^2$:
\begin{equation}
\label{eqn:integral} 
\int_{\Sig} f^*h |\tilde{s}|^2 da = \int_{S^2} h d\mu.
\end{equation}
Apply Lemma \ref{lem:sections} to the line bundle $L^*$ on $S^2$ and Borel measure $\mu$ to find two sections $t_1^*$ and $t_2^*$
satisfying the conclusions of the lemma. Set $\a_1= f^* t_1^*$ and $\a_2= f^* t_2^*$ and consider these sections as maps:
$$\a_i: \G(\xi \otimes E) \to \G(E),$$
via contraction.
Define $s_1 = \a_1(\tilde{s}), s_2 = \a_2(\tilde{s}) \in \G(E)$. Then using $(1)$ of the lemma,
$$|s_1|^2 + |s_2|^2 = ( |\a_1|^2 + |\a_2|^2) |\tilde{s}|^2 \geq |\tilde{s}|^2.$$
Hence,
\begin{equation}
\label{equ:integralbdd}
\int_\Sig |s_1|^2 da + \int_\Sig |s_2|^2 da  \geq \int_\Sig |\tilde{s}|^2 da.
\end{equation}
Using (\ref{eqn:integral}) we have:
$$
\int_{\Sig} |s_2|^2 da  =   \int_{\Sig}  |\a_2|^2 |\tilde{s}|^2 da  =  \int_{S^2} |t_2^*|^2 d\mu
$$
and
$$
\int_{\Sig} |s_1|^2 da  =   \int_{\Sig}  |\a_1|^2 |\tilde{s}|^2 da =   \int_{S^2} |t_1^*|^2 d\mu. 
$$
Therefore using $(3)$ of the lemma,
$$2 \int_{\Sig} |s_1|^2 da \geq \int_{\Sig} |s_2|^2 da.$$
It follows from (\ref{equ:integralbdd}) that,
$$3 \int_\Sig |s_1|^2 da   \geq \int_\Sig |\tilde{s}|^2 da.$$
Using $(2)$ of the lemma as in~\cite{F},
\begin{eqnarray*}
|\bar{\p} s_1|  & = & |\bar{\p} (\a_1(\tilde{s}))|  \\ [.1cm] & =  &|(\bar{\p} \a_1)\tilde{s} + \a_1 (\bar{\p} \tilde{s})|  \\ [.1cm] & = &  |\bar{\p}( f^* t_1^*)\tilde{s} | \\[.1cm]  & \leq & c|\bar{\p} f| |\tilde{s} | \\[.1cm]
& \leq & c \e |\tilde{s}|.
\end{eqnarray*}
Therefore,
\begin{equation}
\int_\Sig |\bar{\p} s_1|^2 da \leq 3(c\e)^2 \int_\Sig | s_1|^2 da < \kappa \int_\Sig | s_1|^2 da.
\end{equation}
by our choice of $\e$. This violates the stability inequality (\ref{equation:stability}), and completes
the proof of Theorem \ref{thm:pic2}.

\begin{rem}
The $\bar{\p}$ operator acting on $s_1$ is determined by the connection on the bundle $\xi^* \otimes \xi \otimes E$
and is therefore not, a priori, equivalent to the $\bar{\p}$ operator determined by the connection on $E$. 
However, the connection on $\xi^* \otimes \xi$ is pulled back from a connection on the  line bundle $L^* \otimes L$ over $S^2$
and is trivial outside a simply connected region $U \subset \Sig$. Recall that, on the trivial line bundle over $S^2$, there is only one holomorphic structure.
Similarly there is a unique  holomorphic structure on $\xi^* \otimes \xi$ that is standard outside $U$.  This
uniqueness implies that the $\bar{\p}$ operator on $\xi^* \otimes \xi$ is equivalent to the standard $\bar{\p}$ operator 
on the trivial line bundle over $\Sig$. Therefore the $\bar{\p}$ operator on $\xi^* \otimes \xi \otimes E$ is equivalent to that
on $E$.
\end{rem}



\bibliographystyle{plain}

\end{document}